\documentclass[preprint,3p,11pt]{elsarticle}

\usepackage{tabls}
\usepackage{cites}
\usepackage{epsf}

\usepackage{hyperref}
\hypersetup{colorlinks=true}
\usepackage[ruled,vlined]{algorithm2e}
\usepackage[caption = false]{subfig}
\usepackage{amsmath}
\usepackage{amssymb}
\usepackage{accents}
\usepackage{xcolor}

\begin{document}
\begin{frontmatter}



\title{Reduced-Memory Methods for Linear Discontinuous Discretization of
the Time-Dependent Boltzmann Transport Equation}

\author[ncsu]{Rylan C. Paye}
\ead{rcpaye@ncsu.edu}

\author[ncsu]{Dmitriy Y. Anistratov}
 \ead{anistratov@ncsu.edu}
\address[ncsu]{Department of Nuclear Engineering,
North Carolina State University, Raleigh, NC}

\author[tamu]{Jim E. Morel}
\address[tamu]{Department of Nuclear Engineering,
Texas A\&M University,  College Station, TX }
 \ead{morel@tamu.edu}
 
\author[lanl]{James  S. Warsa}
\address[lanl]{Los Alamos National Laboratory
  Los Alamos, NM}
 \ead{warsa@lanl.gov}

 % \address[tamu1]{morel@tamu.edu}
 % \address[ncsu1]{rcpaye@ncsu.edu}
%\address[ncsu2]{anistratov@ncsu.edu}
%\address[lanl1]{warsa@lanl.gov}
%---------------------------------------------------------------------------------------------------
\begin{abstract}
In this paper,  new implicit methods with reduced memory are developed for solving the time-dependent Boltzmann transport equation (BTE).
One-group transport problems in 1D slab geometry are considered.
The reduced-memory methods are formulated for the BTE discretized with the linear-discontinuous scheme in space and backward-Euler time
integration method. Numerical results are presented to demonstrate performance of the proposed numerical methods.
%%%%%%%%%%%%%%%%%%%%%%%%%%%%%%%%%
\end{abstract}
 \end{frontmatter}

%---------------------------------------------------------------------------------------------------
\section{Introduction}
 We consider the  time-dependent problems for the one-group linear  Boltzmann transport equation (BTE) in 1D slab geometry
 with isotropic scattering and source   given  by
\begin{equation} \label{bte}
\frac{1}{v}   \partial_t \psi (x, \mu , t)  +
  \mu    \partial_x \psi (x, \mu , t)+\sigma_t(x,t) \psi(x, \mu , t)   =
 \frac{1}{2}\sigma_s(x,t) \int_{-1}^1 \psi(x,\mu',t) d \mu'   + \frac{1}{2}q(x,t) \, ,
\end{equation}
\begin{equation*}
x \in [0,X] \, , \quad  \mu \in [-1,1] \, ,  \quad t \ge 0  \, ,
\end{equation*}
\begin{equation} \label{t-bc}
\psi \big|_{x=0}  =  \psi_{in}^{+} \quad \  \mu>0 \, ,
\quad
\psi \big|_{x=X} = \psi_{in}^{-} \quad \mu<0 \, ,
\quad  \mbox{for} \ t  \ge  0 \, ,\
\end{equation}
\begin{equation} \label{t-ic}
\psi \big|_{t=0} = \psi_{0} \quad \mbox{for}  \quad x \in [0,X] \, , \quad  \mu \in [-1,1]  \, .
\end{equation}
Here $x$ is the spatial position, $\mu$ is the directional cosine of particle motion,
$t$ is time, $\psi$ is the angular flux, $v$ is the particle speed, $q$ is the external source, $\sigma_t$ and $\sigma_s$
are total and scattering cross sections, respectively.

At the end of  each  time step, numerical methods for solving the  time-dependent BTE discretized in phase space and time
requires that the angular flux is stored to   define   the initial condition for the next time step.
This part of the numerical  solution is a set of grid functions of very high dimension that are
specific to the spatial discretization scheme and time integration method.
Various approximation methods have been developed to reduce the associated memory requirements in time-dependent particle transport
problems
\cite{vya-danilova-bnv-1969,pg-dya-jctt,ryan-mc2019,dya-tans-2020,sun-hauck-jsc-2020,dya-jcm-mc2021,einkemmer-jcp-2021,hauck-schnake-arxiv-2022}.
 In this paper,  we  present new approximation methods intended to reduce memory requirements for numerical solution of the time-dependent BTE.
Using  backward-Euler (BE) time integration and the linear discontinuous (LD) spatial discretization,
only the cell-average angular flux is stored between time steps.
In the LD scheme, the slope is represented by the first spatial moment (FSM) of the angular flux on a mesh cell.
At the end of a time step the data on the FSM of the solution is discarded and, in the next time step,
the FSM is reconstructed or approximated in various ways using the cell-average angular flux and/or
the solution to the low-order equations of the second moment (SM) method (used to accelerate iterative convergence)
\cite{lewis-miller-1976,mla-ewl-pne-2002}.
These methods reduce memory allocation for the angular flux by a factor of 2, by computing the FSM
during a time step by calculating the FSM in each mesh cell during every iteration of the high-order transport solve.
In multidimensions this factor is even greater and, furthermore, the approach can have a significant impact
in the context of multiphysics simulations.

The remainder of this paper is organized as follows.
Sec. \ref{disc-BTE} presents the temporal and spatial discretization of the BTE.
Sec. \ref{SM}   briefly describes   the  SM method.
The reduced-memory methods are formulated in Sec. \ref{sec:rmm}.
Sec. \ref{sec:num-res} presents numerical results.
We conclude with a brief discussion in Sec.  \ref{sec:conc}.

%---------------------------------------------------------------------------------------------------
\section{Discretization of   the BTE \label{disc-BTE}}

We define the spatial mesh $\{x_i\}_{i=0}^I$
and angular directions $\{\mu_m\}_{m=1}^M$, where    $i$ is the index of the spatial cell and
   $m$ is the index of the discrete direction.
The LD scheme  in the $i^{th}$ cell,  $x\in [x_i, x_{i+1}]$
for the BTE discretized by the BE time integration method is
defined by the following equations:
 \begin{subequations}  \label{bte-ld-be}
 \begin{equation}  \label{bte-ld0}
  \frac{1}{v\Delta t^n} (\bar \psi_{m,i}^n - \bar\psi_{m,i}^{n-1}) +
\mu_m (\psi_{m,i}^n - \psi_{m,i-1}^n) + \sigma_{t,i}^n\Delta x_i \bar \psi_{m,i}^n =
\frac{\Delta x_i}{2} (\sigma_{s,i}^n \bar \phi_i^n + \bar q_i^n ) \, ,
\end{equation}
\begin{equation}  \label{bte-ld1}
  \frac{1}{v\Delta t^n} (\hat \psi_{m,i}^n - \hat\psi_{m,i}^{n-1}) +
3\mu_m (\psi_{m,i}^n + \psi_{m,i-1}^n - 2 \bar \psi_{m,i}^n) + \sigma_{t,i}^n\Delta x_i \hat \psi_{m,i}^n =
\frac{\Delta x_i}{2}(\sigma_{s,i}^n \hat \phi_i^n + \hat q_i^n) \, ,
\end{equation}
 \begin{equation}\label{bte-ld-aux}
\psi_{m,i-1}^n  =    \bar \psi_{m,i}^n -  \hat \psi_{m,i}^n \,   \  \mbox{for} \ \mu_m < 0 \, , \quad
 \psi_{m,i}^n   =   \bar \psi_{m,i}^n +  \hat \psi_{m,i} ^n\,   \  \mbox{for} \ \mu_m > 0 \, ,
 \end{equation}
  \begin{equation*}
 i \in \mathbb{N}(I), \quad m\in \mathbb{N}(M) \, ,
  \end{equation*}
   \begin{equation}
    \psi_{m,1}^n =   \psi_{in,m}^{+,n}   \  \mbox{for} \ \mu_m > 0 \, , \quad
        \psi_{m,I}^n =   \psi_{in,m}^{-,n}   \  \mbox{for} \ \mu_m < 0 \, ,
   \end{equation}
  \end{subequations}
where
 $n$ is the index of the time layer,
 $\psi_{m,i}^n$ is the  cell-edge angular flux,
\begin{equation}
   \bar \psi_{m,i}^n  = \frac{1}{\Delta x_i}\int_{x_{i-1}}^{x_i} \! \!  \psi_m^{n}dx
\end{equation}
 is the cell-average angular flux,
\begin{equation}
   \hat \psi_{m,i}^n   =    \frac{6}{\Delta x_i^2}  \int_{x_{i-1}}^{x_i}   (x-\bar x_i) \psi_m^{n}dx \,
\end{equation}
is the FSM  of the angular flux,
  $\Delta x_i = x_i -x_{i-1}$, $\bar x_i  =   \frac{1}{2}( x_i +x_{i-1})$,
and $\bar q_{i}^n$ and $\hat q_{i}^n$ are corresponding spatial moments of the source.
 Hereafter we refer to the scheme Eq.~\eqref{bte-ld-be}
 as the LD-BE scheme.
   The discontinuous corner values of the angular flux at the right and left edges of the cell in the LD scheme
are
\begin{equation} \label{psi-L&R}
 \psi_{R,m,i}^n  =  \bar \psi_{m,i}^n +  \hat \psi_{m,i}^n \, ,
\quad
 \psi_{L,m,i}^n  = \bar \psi_{m,i}^n -  \hat \psi_{m,i}^n \, ,
\end{equation}
such that the LD approximation of the solution   in the $i^{th}$ cell defined by the corner values has the form
\begin{equation}
\psi_{m,i}^n(x) = \sum_{\nu=L,R}  \psi_{\nu,m,i}^n B_{\nu,i}(x) \, ,
\end{equation}
with basis function
\begin{equation}
 B_{L,i}(x) =\frac{1}{\Delta x_i} (x_i - x) \, , \quad
 B_{R,i}(x) =\frac{1}{\Delta x_i} (x - x_{i-1}) \, .
\end{equation}

%---------------------------------------------------------------------------------------------------
\section{The Second Moment Method \label{SM}}

To accelerate the iterative solution of the transport problem,
we apply the linear SM method \cite{lewis-miller-1976,mla-ewl-pne-2002}.
The low-order SM (LOSM) equations  are given by
\begin{subequations}  \label{losm}
 \begin{equation} \label{losm0}
  \frac{1}{v} \partial_t \phi (x,t)   +  \partial_x J (x,t)  +
 \sigma_{a}(x,t)  \phi(x,t) = q(x,t) \, ,
\end{equation}
 \begin{equation}\label{losm1}
  \frac{1}{v} \partial_t J (x,t)      +
  \frac{1}{3} \partial_x \phi  (x,t)   + \sigma_t(x,t)  J(x,t)  =
 \partial_x F (x,t)    \, ,
\end{equation}
\end{subequations}
where
 \begin{equation} \label{f-factor}
F =\int_{-1}^1 \Big( \frac{1}{3}   - \mu^2 \Big)  \psi  d \mu \,
\end{equation}
defines the closure for the system of the high-order BTE, Eq.~\eqref{bte}, and the LOSM equations, Eq.~\eqref{losm}.
The LOSM equations are discretized in time with the BE method and spatially discretized
with a scheme that is algebraically consistent with the high-order LD-BE equation, Eq.~\eqref{bte-ld-be}.
The discretized LOSM equations yield the spatial moments of the scalar flux
 \begin{equation}
 \bar \phi_i^n =  \frac{1}{\Delta x_i} \int_{x_{i-1}}^{x_i} \! \!  \phi^{n}dx\, , \quad
 \hat \phi_i^n = \frac{6}{\Delta x_i^2} \int_{x_{i-1}}^{x_i} \! \! (x-\bar x_i) \phi^{n}dx \, ,
 \end{equation}
 the corner values of $\phi_{L,i}^n$ and $\phi_{R,i}^n$, and similar grid functions for the current $J^n$.
The solution algorithm for the two-level SM method on every time step is presented in Algorithm \ref{SM-algorithm}.
Note that algebraic consistency of the discretized high-order and low-order equations is not a requirement for effective
iterative acceleration.
This allows us to introduce various approximations into the LD-BE equation
that reduce memory allocation between time steps and still retain stable and
rapid iterative convergence of the two-level method, even though the discretization
of the high-order and low-order equations are independent as a result of those approximations.

\begin{algorithm}
\DontPrintSemicolon
%\dontprintsemicolon
\For{$n =1,\ldots, N_{t}$}{
$t =  t^n$\;
$s=-1$\;
set $F^{n(1/2)} = F^{n-1}$\;
\nl\While{
$||\phi^{n,s} - \phi^{n,s-1}||_{\infty} < \varepsilon_{\phi}\Big(   \frac{1}{\rho^{s}_{\phi}} - 1 \Big)$
}{\label{InRes1}
$s=s+1$\;
  \For{$s > 1$}
      {
        Solve the high-order transport problem to calculate $\psi^{n,s+1/2}$ using
        $\phi^{n,s}$ to define the right-hand side of the BTE\;
        Compute the SM factor $F^{n,s+1/2}$ using the high-order transport solution $\psi^{n,s+1/2}$\;
      }
      Solve the  LOSM  problem   for $\phi^{n,s+1}$  and $J^{n,s+1}$\;
$\rho^{s}_{\phi} = \frac{||\phi^{n,s} - \phi^{n,s-1}||_{\infty}}
{||\phi^{n,s-1} - \phi^{n,s-2}||_{\infty}}$ \;
}
}
\caption{Iteration Scheme of the Two-Level SM Method \label{SM-algorithm}}
\end{algorithm}

%---------------------------------------------------------------------------------------------------
\section{Reduced-Memory Methods  \label{sec:rmm}}

The reduced-memory methods approximate or reconstruct the FSM
of the angular flux at the previous time step in Eq.~\eqref{bte-ld1}
with
\begin{equation}
\hat \psi_{m,i}^{n-1} \approx \hat \psi_{m,i}^{\star,n-1},
\end{equation}
such that it is modified as follows,
\begin{equation}  \label{bte-ld1-mod}
  \frac{1}{v\Delta t^n} (\hat \psi_{m,i}^n - \hat\psi_{m,i}^{\star,n-1}) +
3\mu_m (\psi_{m,i}^n + \psi_{m,i-1}^n - 2 \bar \psi_{m,i}^n) + \sigma_{t,i}^n\Delta x_i \hat \psi_{m,i}^n =
\frac{\Delta x_i}{2}(\sigma_{s,i}^n \hat \phi_i^n + \hat q_i^n) \, .
\end{equation}
With this modification, the discretization of the LOSM equations and the LD-BE equations,
Eq.~\eqref{bte-ld0}, Eq.~\eqref{bte-ld-aux}, and Eq.~\eqref{bte-ld1-mod}, are not consistent.
Note that the FSM of the scalar flux from the LOSM solution is stored at the end of the time step
and no approximation for the low order scalar flux, $\hat \phi_i^{n-1}$, is introduced in the
LOSM equations.

\subsection{Approximation of the FSM of the Angular Flux}

We consider the following approximations.
\begin{itemize}
 \item The zero-slope approximation. This neglects the slope of the angular flux, that is,
 \begin{equation}
 \hat \psi_{m,i}^{\star,n-1}=0 \, .
\end{equation}

\item The $P_1$ approximation.
  The FSM is based on a $P_1$-expansion of the angular flux using the low-order solution.
  It is defined by
 \begin{equation}
      \hat \psi_{m,i}^{\star \, n-1} =
      \frac{1}{2}\big(\hat\phi_{i}^{n-1}  + 3\mu_m \hat J_{i}^{n-1} \big) \, .
\end{equation}
where $\hat\phi_{i}^{n-1}$ and $\hat J_{i}^{n-1}$  are solutions of the LOSM equations.
\end{itemize}

\subsection{FSM Reconstruction}

To derive a reconstruction for the FSM, we define the approximate corner values in the $i^{th}$ cell
\begin{equation}
 \psi_{L,m,i}^{\star\star \, n-1} = \frac{1}{2}( \bar \psi_{m,i}^{n-1} + \bar \psi_{m,i-1}^{n-1}) \, ,
  \quad   i= 2, \ldots, I \, ,
\end{equation}
\begin{equation}
 \psi_{R,m,i}^{\star\star \, n-1} = \frac{1}{2}( \bar \psi_{m,i+1}^{n-1} + \bar \psi_{m,i}^{n-1}) \, ,
 \quad  i= 1, \ldots, I-1 \,
\end{equation}
using the cell-average values in the two adjacent neighbouring cells on either side of a mesh cell $i$.
The reconstructed FSM is given by
\begin{equation}
\hat \psi_{m,i}^{\star\star \, n-1} =
 \frac{1}{2} (  \psi_{R,m,i}^{\star\star \, n-1} -   \psi_{L,m,i}^{\star\star \, n-1}) \, .
 \quad   i= 1, \ldots, I \, .
\end{equation}
The reconstructed FSM in interior cells is defined by the cell-average values \cite{cockburn-shu-1989}:
\begin{equation} \label{slope-interior}
\hat \psi_{m,i}^{\star\star \, n-1} =
 \frac{1}{4} ( \bar \psi_{m,i+1}^{n-1} - \bar\psi_{m,i-1}^{n-1}) \, ,
  \quad   i= 2, \ldots, I-1 \, .
\end{equation}
In the boundary cells, we define the corner values as
\begin{equation}
\ \psi_{L,m,1}^{\star\star \, n-1} =
\begin{cases}
 \frac{1}{2}( \bar \psi_{m,1}^{n-1}+  \psi_{in,m}^{+,n-1}) \, , &
  \mu >0 \,  ,\\
  \bar \psi_{m,1}^{n-1}\, , &
  \mu < 0 \, ,
  \end{cases}
\end{equation}
\begin{equation}
 \psi_{R,m,I}^{\star\star \, n-1} =
\begin{cases}
 \bar \psi_{m,I}^{n-1} \, , &
  \mu >0 \, , \\
  \frac{1}{2}( \bar \psi_{m,I}^{n-1} +   \psi_{in,m}^{-,n-1}) \, , &
  \mu < 0 \, .
  \end{cases}
\end{equation}
This gives rise to
\begin{subequations}\label{sr&bcells-corner-v}
\begin{equation}
\hat \psi_{m,1}^{\star\star,n-1} =
\begin{cases}
 \frac{1}{4}( \bar \psi_{m,2}^{n-1} -  \psi_{in,m}^{+,n-1}) \, , &
  \mu >0 \,  ,\\
  \frac{1}{4}(\bar \psi_{m,2}^{n-1} - \bar \psi_{m,1}^{n-1}) \, , &
  \mu < 0 \, ,
  \end{cases}
\end{equation}
\begin{equation}
\hat\psi_{m,I}^{\star\star,n-1} =
\begin{cases}
 \frac{1}{4}(\bar \psi_{m,I}^{n-1}  - \bar \psi_{m,I-1}^{n-1} )\, , &
  \mu >0 \, , \\
  \frac{1}{4}(\psi_{in,m}^{-,n-1} -  \bar \psi_{m,I-1}^{n-1}  ) \, , &
  \mu < 0 \, .
  \end{cases}
\end{equation}
\end{subequations}

To obtain the FSM reconstruction for $\hat \psi_{m,i}^{\star \, n-1}$ in Eq.~\eqref{bte-ld1-mod},
we apply a slope limiter (SL) to $\hat\psi_{m,i}^{\star\star,n-1}$ given by
\cite{cockburn-shu-1989,ryan-rob-jcp-2008}
 \begin{equation} \label{minmod-sl}
 \hat \psi_{m,i}^{\star \, n-1} =  minmod \Big(\hat \psi_{m,i}^{\star\star \ , n-1},
 \frac{1}{2} \big( \bar \psi_{m,i}^{ n-1} - \bar \psi_{m,i-1}^{ n-1} \big),
 \frac{1}{2} \big( \bar \psi_{m,i+1}^{ n-1} - \bar \psi_{m,i}^{ n-1} \big)
 \Big) \, ,
 \end{equation}
 where
  \begin{equation}
  minmod  (f_1,f_2,f_3) =
  \begin{cases}
  sign(f_1)\min(|f_1|,|f_2|,|f_3|) & \mbox{if} \
   sign(f_1) =   sign(f_2)=  sign(f_3) \, ,\\
  0  & \mbox{otherwise}  \, .
    \end{cases}
\end{equation}
 Hereafter we refer to this as the slope reconstruction (SR) with slope limiting (SL), or SR-SL, scheme.

\subsection{Approximation of the Rate of Change of the FSM of the Angular Flux}

Another group of methods defines
$\hat \psi_{m,i}^{\star , n-1}$
by means of a factor $\beta^n$ that approximates the rate of change of the solution over the $n^{th}$ time step.
The factor $\beta^n$ is formulated by applying various grid functions of the low-order solutions.
\begin{itemize}
\item The $\beta_{\bar \phi}\,$-approximation. This is formulated in terms of the change in the cell-average scalar flux
  from the LOSM equations and defined by
 \begin{equation}
  \beta_{i}^n = \frac{\bar\phi_{i}^{n-1}}{\bar \phi_{i}^n} \,,
 \end{equation}
and the approximate FSM of the angular flux has the form
 \begin{equation} \label{beta-psi}
 \hat \psi_{m,i}^{\star , n-1} = \beta_{i}^n   \hat \psi_{m,i}^n \, .
  \end{equation}
The approximation Eq.~\eqref{beta-psi} is introduced in Eq.~\eqref{bte-ld1-mod} to get
\begin{equation}  \label{bte-ld1-mod-2}
  \frac{1}{v\Delta t^n} (1-  \beta_{i}^n)\hat \psi_{m,i}^n +
3\mu_m (\psi_{m,i}^n + \psi_{m,i-1}^n - 2 \bar \psi_{m,i}^n) + \sigma_{t,i}^n\Delta x_i \hat \psi_{m,i}^n =
\frac{\Delta x_i}{2}(\sigma_{s,i}^n \hat \phi_i^n + \hat q_i^n) \, .
\end{equation}

\item The $\beta_{LR}$-approximation. This is based on the change in the corner values of the scalar flux from the
  LOSM equations and defined by
\begin{subequations}  \label{beta-LR-app}
      \begin{equation} \label{psi-beta-LR}
\hat \psi_{m,i}^{\star \, n-1} = \frac{1}{2} (\beta_{R,i}^n \psi_{R,m,i}^{ n}  - \beta_{L,i}^n\psi_{L,m,i}^{n} ) \, ,
\end{equation}
 where
      \begin{equation} \label{beta-LR}
    \beta_{R,i}^n = \frac{\phi_{R,i}^{n-1}}{\phi_{R,i}^n} \, , \quad
        \beta_{L,i}^n = \frac{\phi_{L,i}^{n-1}}{\phi_{L,i}^n} \, .
\end{equation}
\end{subequations}
  \end{itemize}
The factors $\beta_{i}^n$ and $\beta_{\gamma,i}^n$ ($\gamma=L,R$)
couple the high-order and low-order equations nonlinearly,
and the two-level SM method in Algorithm \ref{SM-algorithm} becomes a nonlinear iteration.

%---------------------------------------------------------------------------------------------------
\section{Numerical Results \label{sec:num-res}}

Numerical results comparing the various reduced-memory methods are shown for two test problems.
The reference solution is that of the two-level SM method without any memory-reduction approximations
made in the discretized LD-BE scheme in Eq.~\eqref{bte-ld-be}.

\subsection{Test A}

\begin{figure}[t!]
	\centering \hspace*{-.5cm}
	\subfloat[\label{test-a-phi-bar-sol}  $\phi(x,t^n)$.]{\includegraphics[width=.525\textwidth]{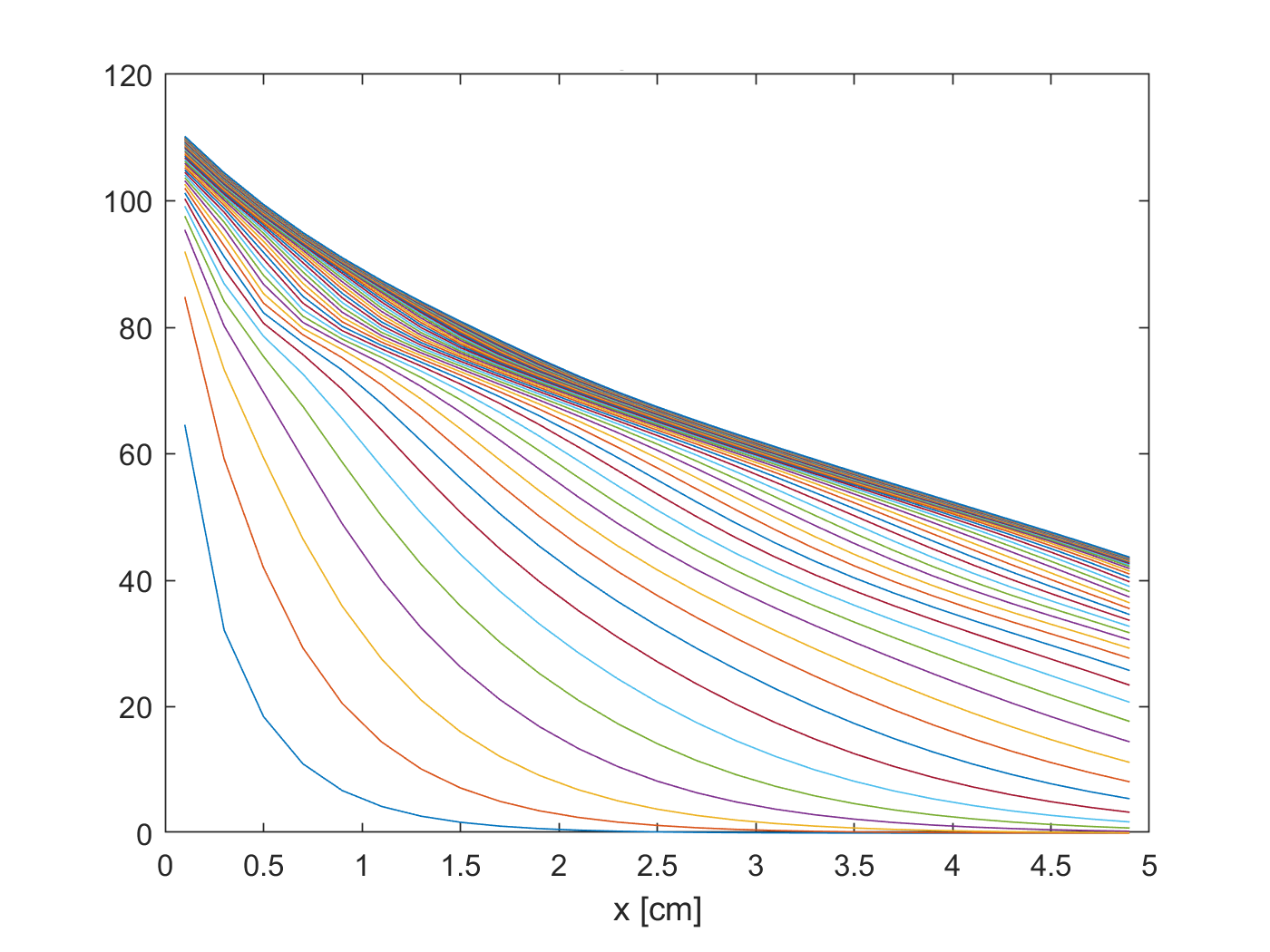}}
	\subfloat[\label{test-a-phi-bar-lambda}$\lambda =\frac{1}{\phi}\partial_t \phi \big|_{t^n}$.]
{\includegraphics[width=.525\textwidth]{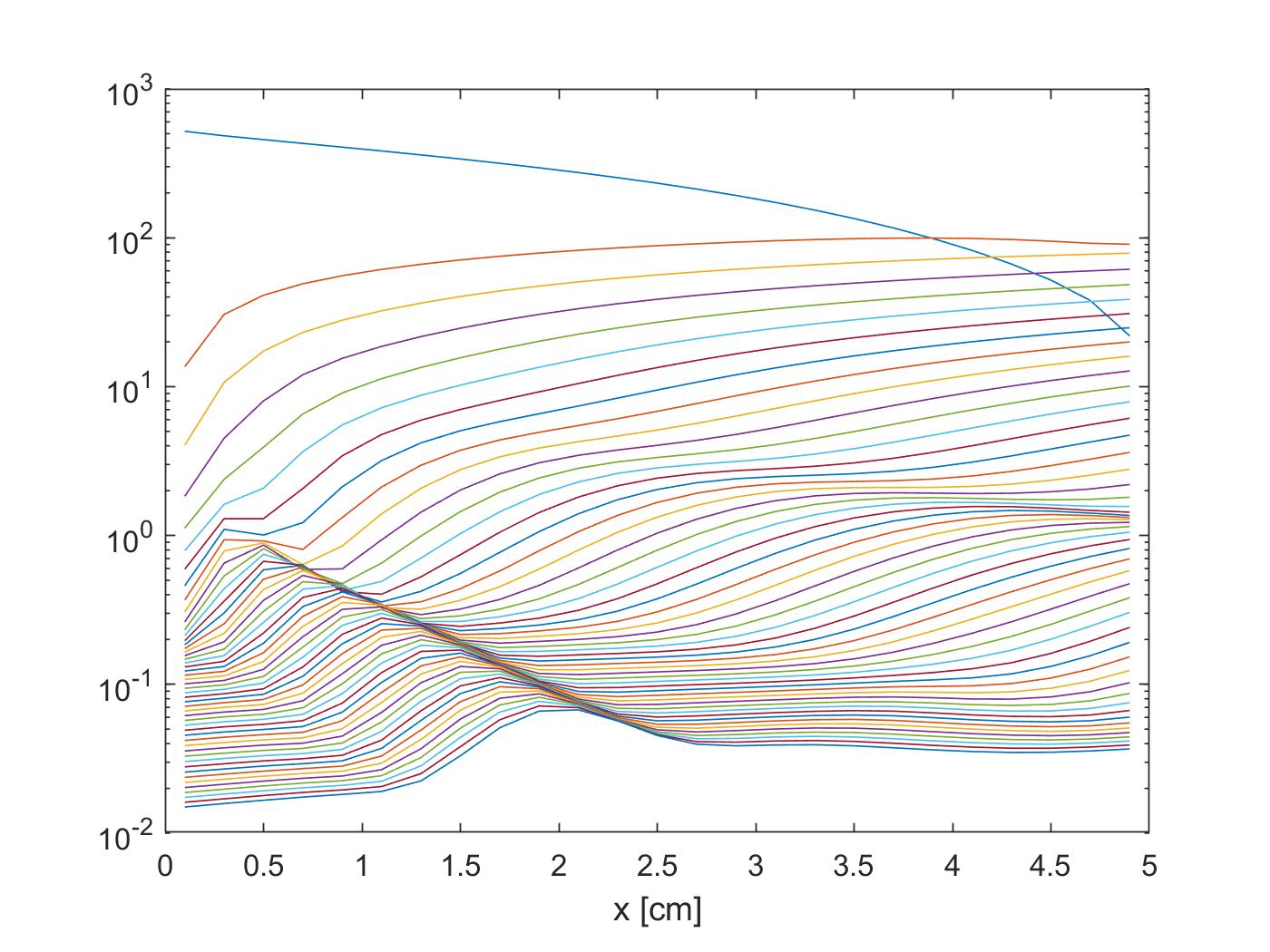}}
	\caption{ \label{test-a-sol} Test A. Numerical solution and its  relative change rate.}
%\end{figure}

%\begin{figure}[t!]
	\centering \hspace*{-.5cm}
	\subfloat[\label{test-a-phi-bar} Relative error in $\bar \phi$.]{\includegraphics[width=.525\textwidth]{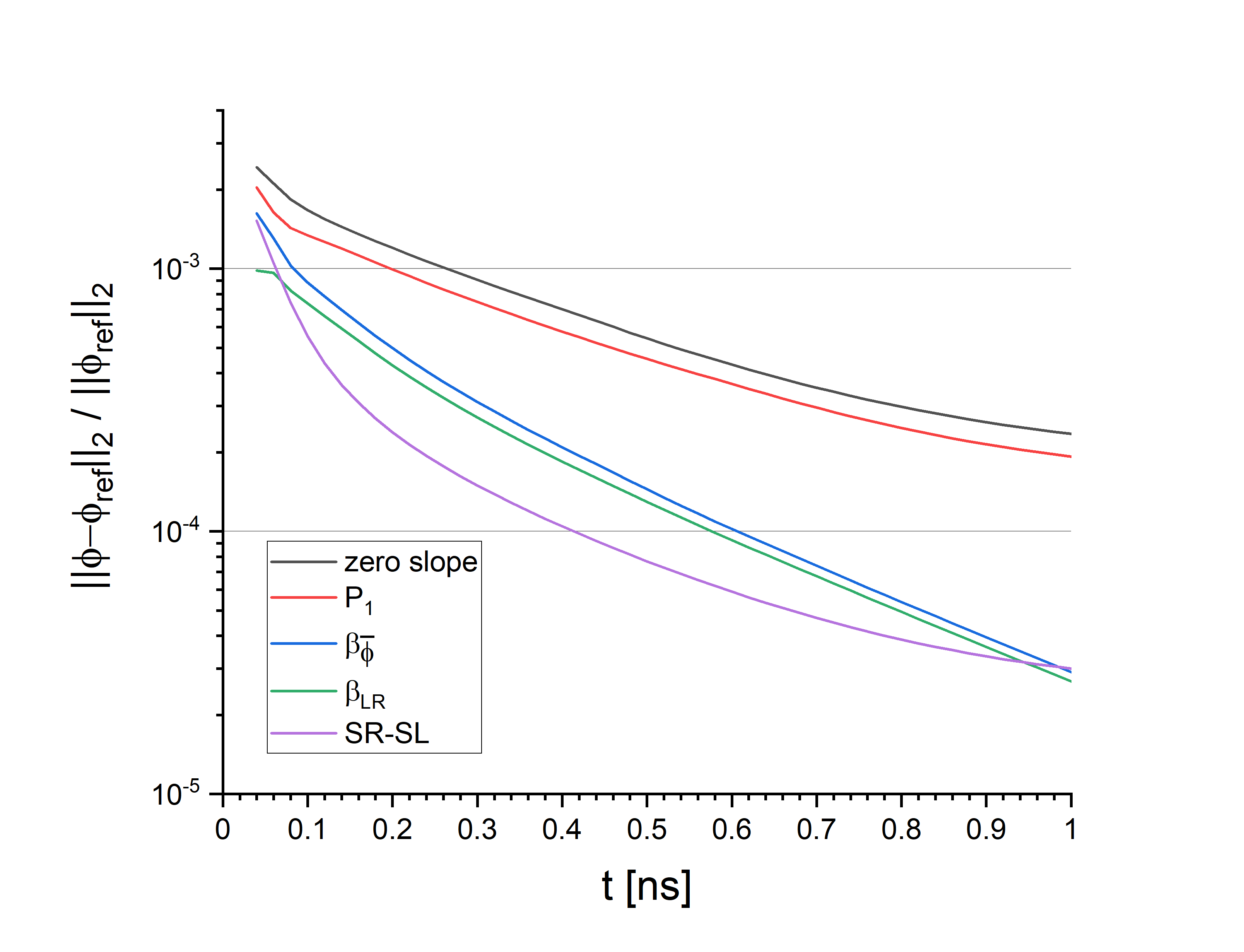}}
	\subfloat[\label{test-a-phi-hat} Relative error in $\hat \phi$.]{\includegraphics[width=.525\textwidth]{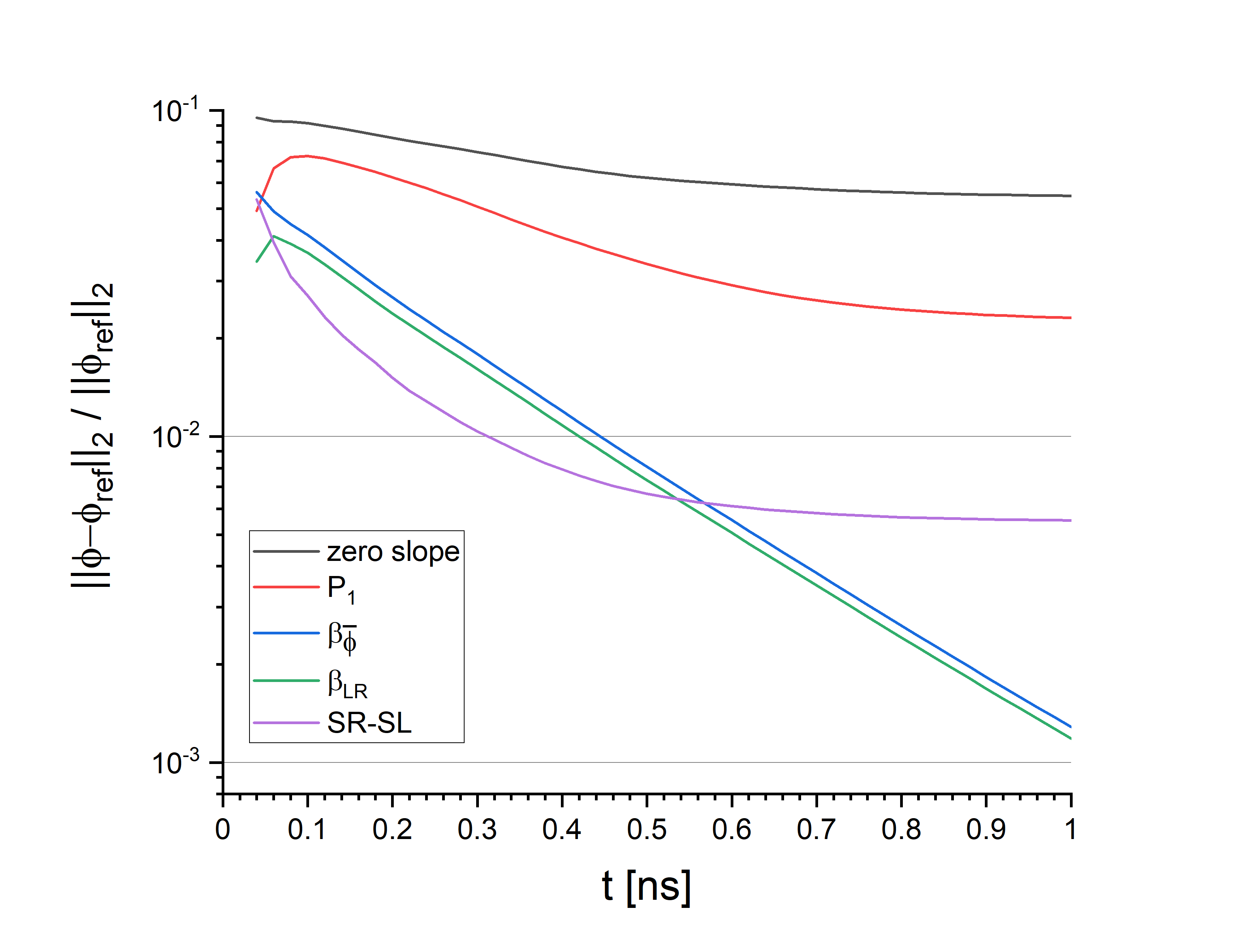}}
	\caption{ \label{test-a} Test A. Relative error in LOSM solutions with reduced-memory methods.}
\end{figure}

This problem models typical high-energy photon transport
with high change rates in the solution.
The test  is defined on the spatial domain $x \in [0,5]\text{cm}$
with $\sigma_t=10^{-1}\text{cm}^{-1}$,
 $\sigma_s  =\ 5 \times 10^{-2}\text{cm}^{-1}$,
$q  =  0$,
$\psi_{in}^+=10^2$, $\psi_{in}^-=0$,
$\psi_0=10^{-3}$,
$v=c= 3\times 10^{10}\frac{\text{cm}}{\text{s}}$,
for $t \in  [0,1]\text{ns}$.
The convergence criterion $\varepsilon_{\phi} =10^{-8}$.
$\psi$ is the total radiation intensity with units
$\big[2\pi \times 10^{21} \frac{\text{erg}}{\text{cm}^2 \cdot \text{s} \cdot \text{ster}}\big]$.
The spatial mesh comprises 100 equal-width intervals.
The quadrature set is a double $S_4$ Gauss-Legendre so that the number of directions is $M=8$.
The time step is constant {$\Delta t=2\times10^{-2}$ ns}.
Figure~\ref{test-a-sol} shows the scalar flux  and the estimated rate of change $\lambda = \frac{1}{\phi} \partial_t \phi$  on every time step in this problem.

Figure~\ref{test-a} presents the relative difference in the numerical solution of the approximate methods
in the 2-norm compared to the reference numerical solution of the LD-BE scheme on the given phase-space grid.
The error in $\bar \phi$ and $\hat \phi$ are shown in Figures \ref{test-a-phi-bar} and \ref{test-a-phi-hat}, respectively.
Note that all the methods use  the  initial conditions given in the problem, i.e. $\hat \psi_{m,i}^{0}$, on the first time step. Thus the error
equals zero at the first time step.

The zero-slope approximation indicates the worst-case for comparison of the other methods.
The $P_1$ approximation of $\hat \psi^{n-1}$ demonstrates  limitation of this natural technique.
The $\beta_{\bar \phi}\,$-approximation shows improvement compared to the  $P_1$ approximation.
The solution of the method with the
$\beta_{LR}$-approximation is more accurate than the one with
the $\beta_{\bar \phi}$-approximation.
The SR-SL scheme generates the most accurate the cell-average scalar flux $\bar \phi$.
The errors in the FSM scalar flux, $\hat \phi$, are larger that those in  $\bar \phi$.
We note that the methods with  $\beta_{\bar \phi}$- and $\beta_{LR}$-approximations compete
with the SR-SL scheme in calculations of $\hat \phi$.
In this test, the number of iterations is the same on each time step.
Table ~\ref{test-a-iter} shows the number of iterations for each method.

\begin{table}[h]
	\centering
	\caption{\label{test-a-iter} Test A. Number of iterations (the same in every time step).}
\medskip
	\begin{tabular}{|l||c|c|c|c|c|c|}
 		\hline
Method & ref  & zero slope  & $P_1$ & SR-SL &  $\beta_{\bar \phi}$	 	&$ \beta_{LR}$	\\ \hline
Iterations & 5	&	5	&	5	&	5	&	5	&	8	 \\ \hline
\end{tabular}
\end{table}

\subsection{Test B}

\begin{figure}[h!]
	\centering \hspace*{-.5cm}
	\subfloat[\label{test-b-phi-bar-sol}  $\phi(x,t^n)$.]{\includegraphics[width=.525\textwidth]{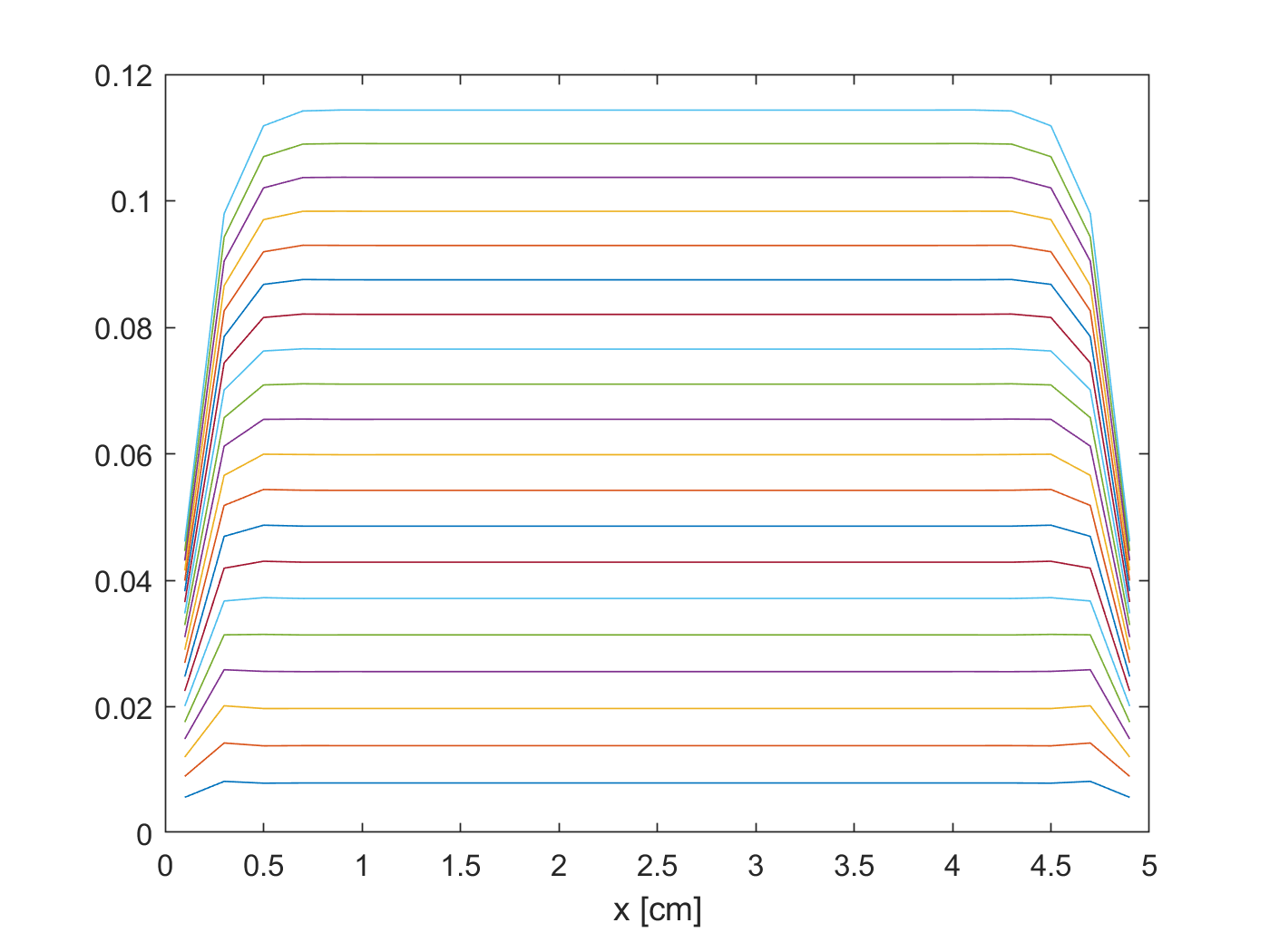}}
	\subfloat[\label{test-b-phi-bar-lambda}$\lambda =\frac{1}{\phi}\partial_t \phi \big|_{t^n}$.]{\includegraphics[width=.525\textwidth]{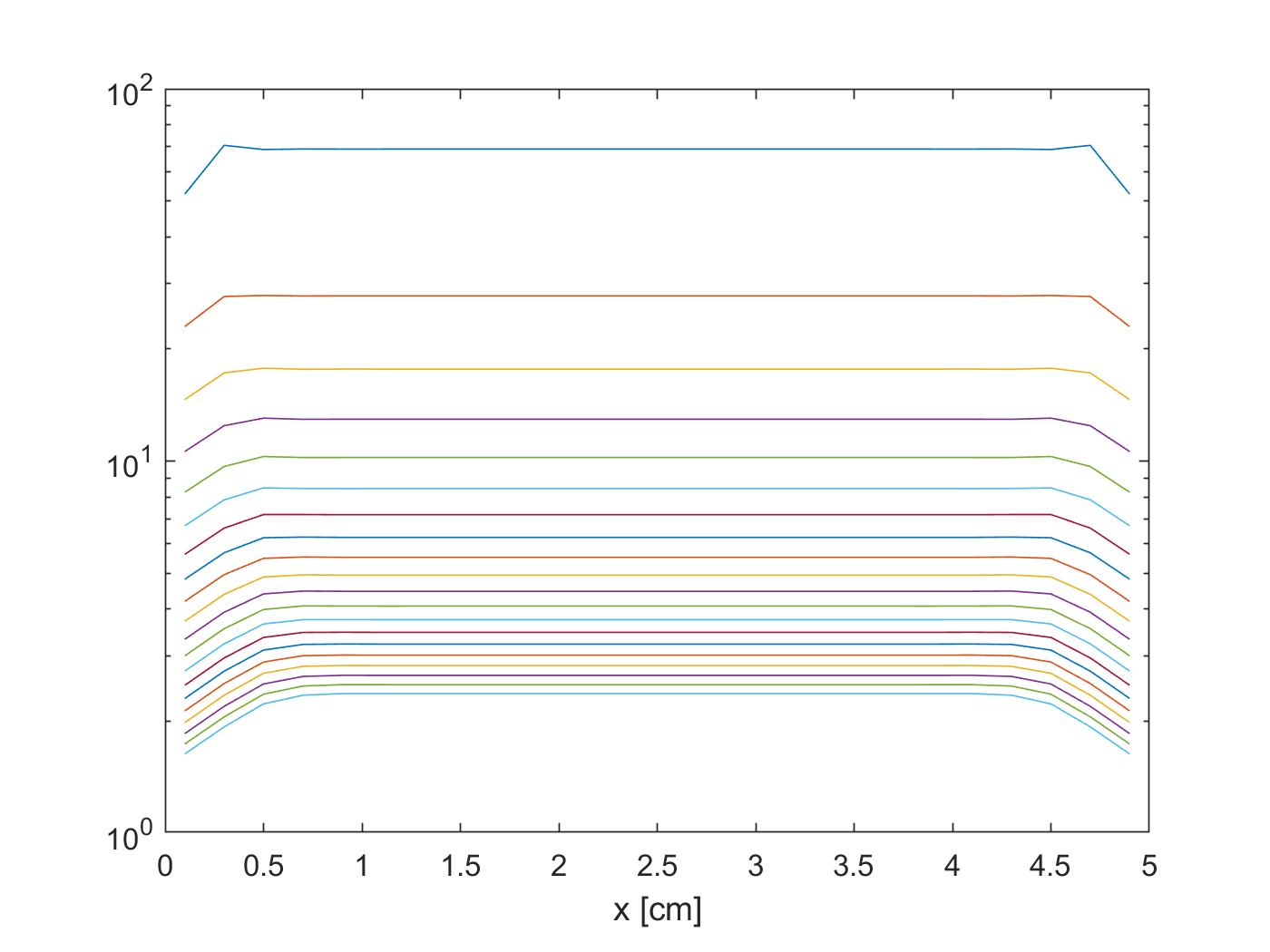}}
	\caption{ \label{test-b-sol} Test B. Numerical solution and its  relative change rate.}
\end{figure}

This is a highly diffusive problem, defined on the spatial domain $x \in [0,5]\text{cm}$ with
$\sigma_t = 10^{2}\text{cm}^{-1}$,
$\sigma_a =10^{-2} \text{cm}^{-1}$,
$q=10^{-2}$,
$\psi_{in}^+=0$, $\psi_{in}^-=0$, $\psi_0=10^{-3}$,
$v=c=3\times10^{10}\frac{\text{cm}}{\text{s}}$,
$t \in [0,0.4] \text{ns}$.
The convergence criterion $\varepsilon_{\phi} =10^{-8}$.
The spatial mesh is uniform with 25 intervals.
The double $S_4$ Gauss-Legendre quadrature set is used and hence $M=8$.
The time step is constant and {$\Delta t=2\times10^{-2}$ ns}.
Figure~\ref{test-b-sol} shows the scalar flux  and the estimated rate of change  on every time step in this test.

\begin{figure}[h!]
	\centering \hspace*{-.5cm}
	\subfloat[\label{test-b-phi-bar} Relative error in $\bar \phi$.]{\includegraphics[width=.525\textwidth]{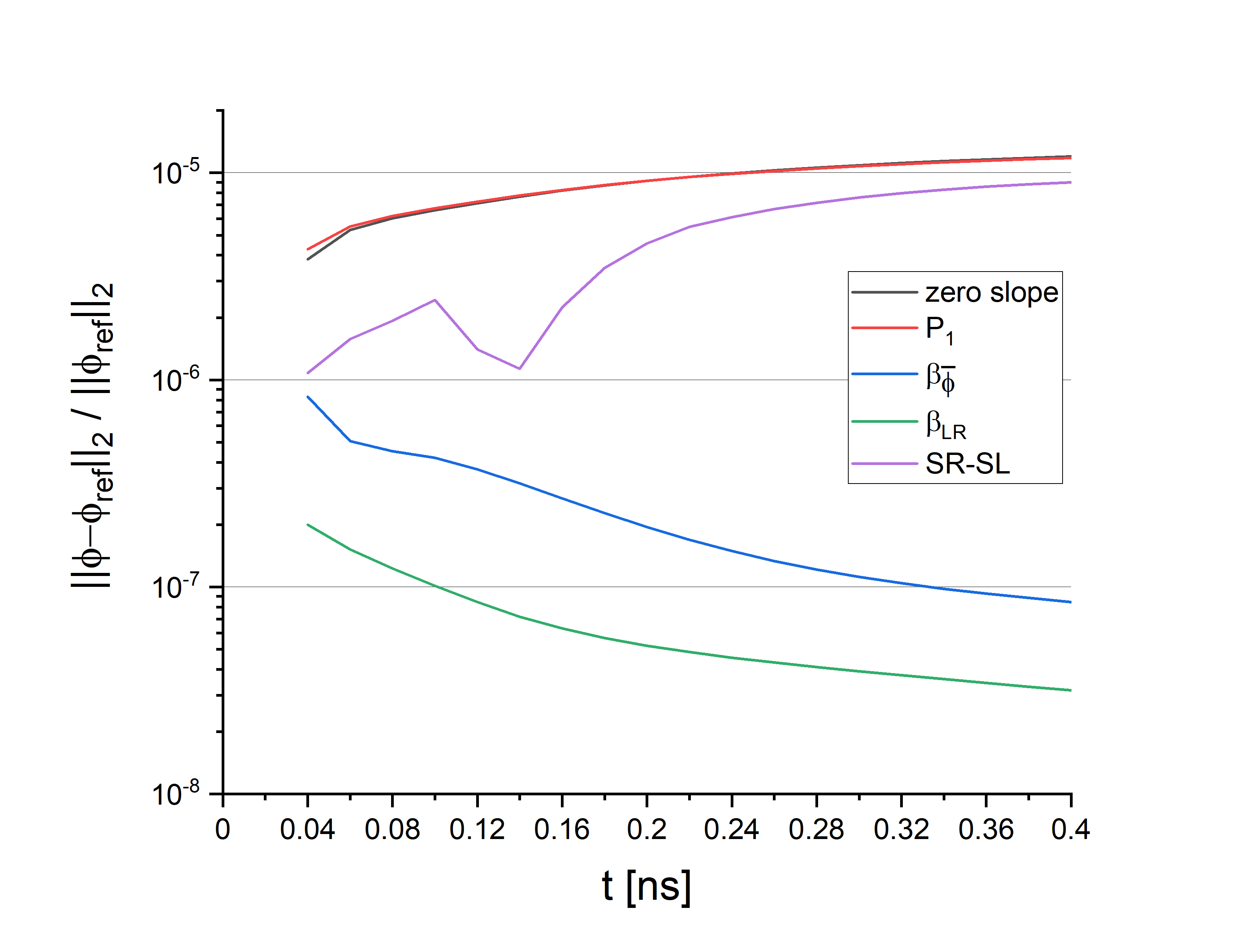}}
	\subfloat[\label{test-b-phi-hat} Relative error in  $\hat \phi$.]{\includegraphics[width=.525\textwidth]{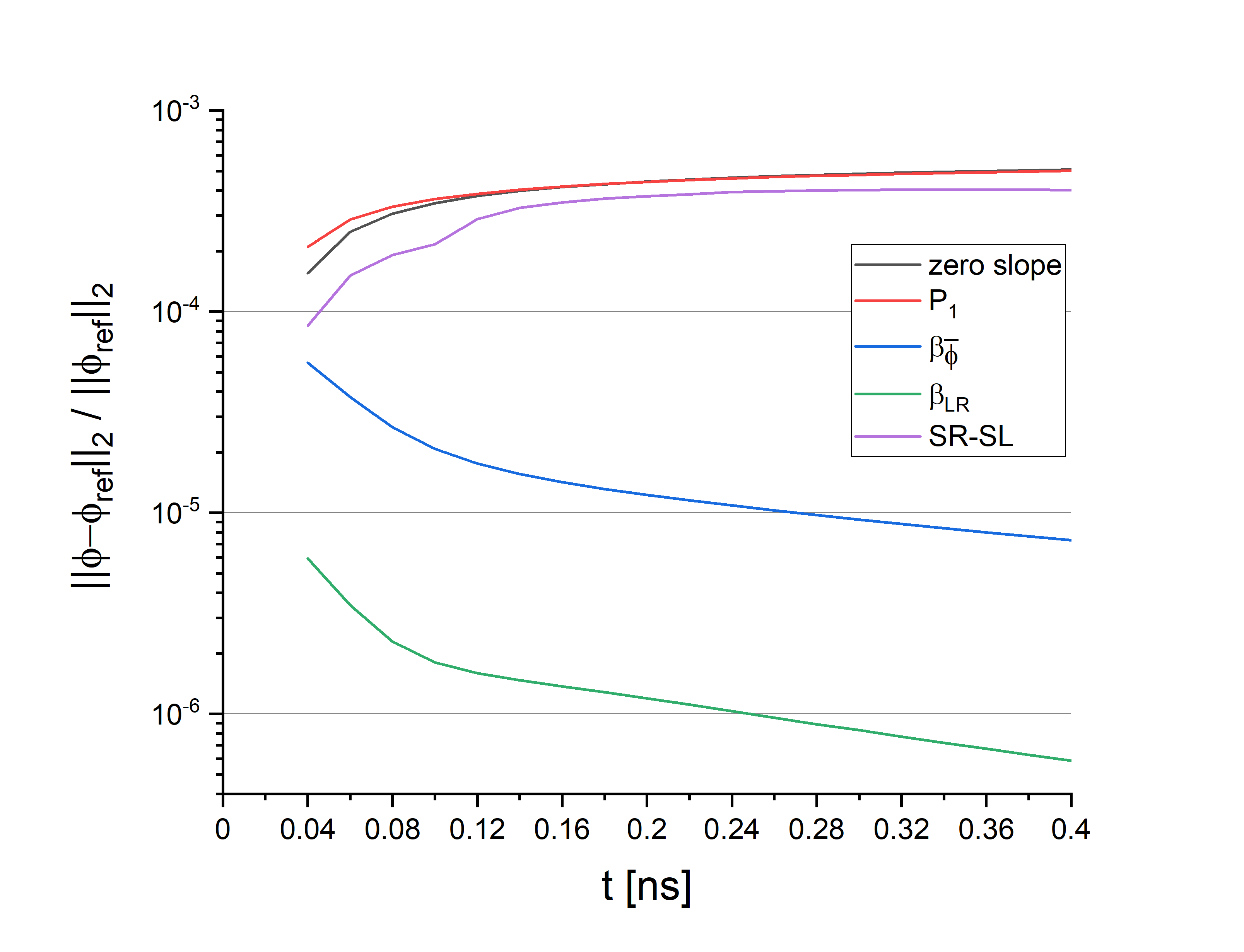}}
	\caption{ \label{test-b} Test B. Relative error in LOSM solutions with reduced-memory methods.}
\end{figure}

\begin{figure}[h!]
 \begin{center}
 \includegraphics[scale=0.45]{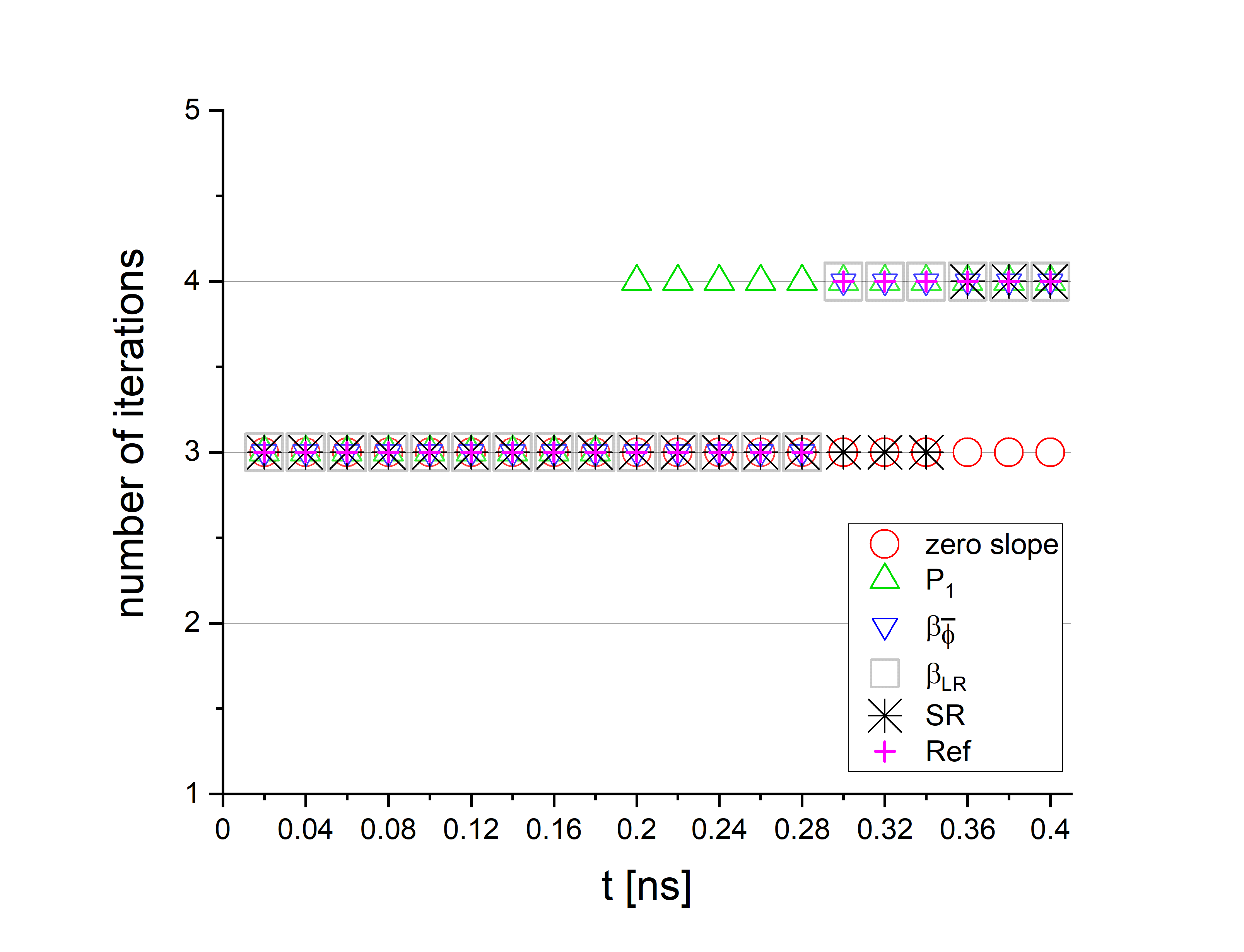}
\caption{\label{test-e-iteration}  Test B. Number of iterations per time step.}
\end{center}
 \end{figure}

Figures \ref{test-b-phi-bar} and \ref{test-b-phi-hat} present errors in $\bar \phi$ and $\hat \phi$, respectively.
The results show that all approximation methods preserve the asymptotic thick-diffusion limit.
The schemes with  $\beta_{\bar \phi}$- and $\beta_{LR}$-approximations have smallest errors.
The number of transport iterations versus time step are shown in Figure \ref{test-e-iteration}.

%---------------------------------------------------------------------------------------------------
\section{Conclusions \label{sec:conc}}

We developed a family of approximation methods that aim to reduce memory resource allocation for
numerical solution of the time-dependent BTE discretized with the LD scheme in space and BE in time.
Numerical results were shown that indicate the SR-SL method and the schemes with $\beta$-approximations
are promising for use in multi-dimensional geometries and nonlinear thermal radiative transfer problems.
The $\beta_{\bar \phi}$- and $\beta_{LR}$-approximations of $\hat \psi_{m,i}^{(n-1)}$ by means of solution
of the LOSM equations lead to nonlinear coupling between the discretized BTE and LOSM equations.
As a result, the transport iterations involving the high-order and low-order equations become a
nonlinear iteration scheme that needs further analysis.
The LD scheme, not unexpectedly, does not preserve positivity of the transport solution.
In general, non-positivity of the numerical transport solution will affect the solution of the formulated
approximation methods and will be addressed in future research.

\section*{Acknowledgements}
Los Alamos Report LA-UR-23-25062. The work of the second and third authors (DYA and JEM)  was funded by the Joint Center for Resilient National  Security  of  Los Alamos National Laboratory and Texas A\&M University System.
The work of the fourth author (JSW) was supported by the U.S. Department of Energy through the Los Alamos National Laboratory. Los Alamos National Laboratory is operated by Triad National Security, LLC, for the National Nuclear Security Administration of U.S. Department of Energy (Contract No. 89233218CNA000001).

\bibliographystyle{elsarticle-num}
\bibliography{rp-dya-jm_jw_RMM-for-BTE}

\end{document}